\documentclass[11pt,leqno,a4paper]{article}

\usepackage{amsmath, amssymb, amsthm, amsfonts}
\usepackage{mathrsfs}

\theoremstyle{plain}
\newtheorem{thm}{Theorem}

\newtheorem{lem}{Lemma}

\theoremstyle{remark}

\numberwithin{equation}{section}

\begin{document}

\title{Approximation and bounds for the Wallis ratio}
\author{ Xu You }

\maketitle

\footnote[0]{2010 Mathematics Subject Classification: 33B15, 26A48, 26D07}
\footnote[0]{Key words and phrases: Wallis ratio, Gamma function, Inequalities, Multiple-correction method}

\begin{abstract}
In this paper, we present an improved continued fraction approximation of the Wallis ratio. This approximation is fast in comparison with the recently discovered asymptotic series. We also establish the double-side inequality related to this approximation. Finally, some numerical computations are provided for demonstrating the superiority of our approximation.

\end{abstract}

\section{Introduction}

The Wallis ratio is defined as
$$W(n)=\frac {(2n-1)!!}{(2n)!!}=\frac 1 {\sqrt \pi}\frac {\Gamma (n+\frac 12)}{\Gamma (n+1)},$$
where $\Gamma$ is the classical Euler gamma function which may be defined by
$$\Gamma(x)=\int _0 ^\infty {t^{x-1}e^{-t}}dt, \textmd{Re}(x)>0.$$
The study and applications of $W(n)$ have a long history, a large amount of literature, and a lot of new results. For detailed information, please refer these papers \cite {chen1, chen2, qi1, qi2} and references cited therein.

Chen and Qi \cite{chen1} presented the following inequalities for the Wallis ratio for every natural number $n$:
\begin{align}
\frac 1 {\sqrt {\pi(n+\frac 4\pi -1)}}\leq \frac{(2n-1)!!}{(2n)!!}\leq \frac 1{\sqrt {\pi (n+\frac 14)}},
\end{align}
where the constants $\frac 4\pi -1$ and $\frac 14$ are the best possible.

Guo, Xu and Qi proved in \cite {guo1} that the double inequality
\begin{align}
\sqrt {\frac e \pi}\left(1-\frac 1{2n}\right)^n\frac {\sqrt {n-1}}n<W(n)\leq\frac 43\left(1-\frac 1{2n}\right)^n \frac {\sqrt {n-1}}n,
\end{align}
for $n\geq2$ is valid and sharp in the sense that the constants $\sqrt {\frac e \pi}$ and $\frac 43$ are best possible. They also proposed the approximation formula
\begin{align}
W(n)\sim \chi (n):=\sqrt {\frac e \pi}\left(1-\frac 1{2n}\right)^n\frac {\sqrt {n-1}}n, n\rightarrow\infty.
\end{align}

Recently, Qi and Mortici \cite{qi3} improved the approximation formula (1.3) as following,
\begin{align}
W(n)\sim \sqrt {\frac e \pi}\left(1-\frac 1{2\left(n+\frac 13\right)}\right)^{n+\frac 13}\frac 1{\sqrt {n}}.
\end{align}

Motivated by these works, in this paper we will apply the \emph{multiple-correction method} \cite {cao2, cao3, cao4} to construct an improved continued fraction asymptotic expansion for the Wallis ratio as follows:

\begin{thm} For the Wallis ratio $W(n)=\frac {(2n-1)!!}{(2n)!!}$, we have
\begin{align}
W(n)\sim \sqrt \frac e \pi \left(1-\frac 1{2\left(n+\frac 13\right)}\right)^{n+\frac 13}\frac 1{\sqrt {n}}\exp \left(\frac 1{n^2}\frac {a_1}{n+b_1+\frac {a_2}{n+b_2+\frac {a_3}{n+b_3}}}\right),
\end{align}
where $a_1=\frac 1 {144}, b_1=\frac 1{60}; a_2=\frac {781}{3600}, b_2=-\frac{4309}{109340}; a_3=\frac {51396085}{89664267}, b_3=\frac {25682346121}{449571834712}.$
\end{thm}

Using Theorem 1, we provide some inequalities for the Wallis ratio.

\begin{thm} For every integer $n>1$, it holds:
\begin{align}
&\sqrt \frac e \pi \left(1-\frac 1{2\left(n+\frac 13\right)}\right)^{n+\frac 13}\frac 1{\sqrt {n}}\exp \left(\frac 1{n^2}\frac {a_1}{n+b_1+\frac {a_2}{n+b_2+\frac {a_3}{n+b_3}}}\right)\\\nonumber
&>W(n)>
\sqrt \frac e \pi \left(1-\frac 1{2\left(n+\frac 13\right)}\right)^{n+\frac 13}\frac 1{\sqrt {n}}\exp \left(\frac 1{n^2}\frac {a_1}{n+b_1+\frac {a_2}{n+b_2}}\right),
\end{align}
where $a_1=\frac 1 {144}, b_1=\frac 1{60}; a_2=\frac {781}{3600}, b_2=-\frac{4309}{109340}; a_3=\frac {51396085}{89664267}, b_3=\frac {25682346121}{449571834712}.$
\end{thm}

To obtain Theorem 1, we need the following lemma which was used in \cite{mor1, mor2, mor3} and is very useful for constructing asymptotic expansions.
\begin{lem}
If the sequence $(x_n)_{n\in \mathbb{N}}$ is convergent to zero and there
exists the limit
\begin{align}
\lim_{n\rightarrow +\infty}n^s(x_n-x_{n+1})=l\in [-\infty,+\infty]
\end{align}
with $s>1$, then
\begin{align}
\lim_{n\rightarrow +\infty}n^{s-1}x_n=\frac{l}{s-1}.
\end{align}
\end{lem}

Lemma 1 was proved by Mortici in \cite {mor1}. From Lemma 1, we can see that the speed of convergence of the sequences $(x_n)_{n\in \mathbb{N}}$ increases together with the values $s$ satisfying (1.8).

The rest of this paper is arranged as follows. In section 2, we will apply the \emph{multiple-correction method} to construct a new asymptotic expansion for the Wallis ratio and prove Theorem 1 by the \emph{multiple-correction method}. In section 3, we established the double-side inequality for the Wallis ratio. In section 4, we give some numerical computations which demonstrate the superiority of our new series over some formulas found recently.

\section{Proof of Theorem 1}
According to the argument of Theorem 5.1 in \cite {qi3} , we can introduce a sequence $(u(n))_{n\geq1}$ by the relation
\begin{align}
W(n)=\sqrt \frac e \pi \left(1-\frac 1{2\left(n+\frac 13\right)}\right)^{n+\frac 13}\frac 1{\sqrt {n}}\exp u(n),
\end{align}
and to say that an approximation $W(n)\sim\sqrt \frac e \pi \left(1-\frac 1{2\left(n+\frac 13\right)}\right)^{n+\frac 13}\frac 1{\sqrt {n}}$ is better if the speed of convergence of $u(n)$ is higher.

\noindent{\bf (Step 1) The initial-correction.}
When $n\rightarrow\infty$, we define a sequence $(u_0(n))_{n\geq1}$
\begin{align}
W(n)=\sqrt \frac e \pi \left(1-\frac 1{2\left(n+\frac 13\right)}\right)^{n+\frac 13}\frac 1{\sqrt {n}}\exp u_0(n).
\end{align}
From (2.2), we have
\begin{align}
u_0(n)=\ln W(n)-\ln \sqrt \frac e \pi- \left(n+\frac 13\right)\ln \left(1-\frac 1{2\left(n+\frac 13\right)}\right)-\frac 12\ln\frac 1n.
\end{align}
Thus,
\begin{align}
u_0(n)-u_0(n+1)=&\ln \frac {2n+2}{2n+1}- \left(n+\frac 13\right)\ln \left(1-\frac 1{2\left(n+\frac 13\right)}\right)-\frac 12\ln \frac 1n\\\nonumber
&+\left(n+\frac 43\right)\ln \left(1-\frac 1{2\left(n+\frac 43\right)}\right)+\frac 12\ln \frac 1{n+1}.
\end{align}
Developing (2.4) into power series expansion in $1/n$, we have
\begin{align}
u_0(n)-u_0(n+1)=\frac 1{48} \frac 1 {n^4} + O(\frac 1{n^5}),
\end{align}
By Lemma 1, we know that the rate of convergence of the sequence $(u_0(n))_{n\geq1}$ is $n^{-3}$.

\noindent{\bf (Step 2) The first-correction.} We define the sequence $(u_1(n))_{n\geq1}$ by the relation
\begin{align}
W(n)=\sqrt \frac e \pi \left(1-\frac 1{2\left(n+\frac 13\right)}\right)^{n+\frac 13}\frac 1{\sqrt {n}}\exp \left(\frac 1{n^2}\frac {a_1}{n+b_1}\right)\exp u_1(n).
\end{align}
From (2.6), we have
\begin{align}
u_1(n)-u_1(n+1)=&\ln \frac {2n+2}{2n+1}- \left(n+\frac 13\right)\ln \left(1-\frac 1{2\left(n+\frac 13\right)}\right)-\frac 12\ln \frac 1n +\frac 12\ln \frac 1{n+1}\\\nonumber
&+\left(n+\frac 43\right)\ln \left(1-\frac 1{2\left(n+\frac 43\right)}\right)-\frac 1{n^2}\frac {a_1}{n+b_1}+\frac 1{(n+1)^2}\frac {a_1}{n+1+b_1}.
\end{align}
Developing (2.7) into power series expansion in $1/n$, we have
\begin{align}
u_1(n)-u_1(n+1)=\left(\frac 1 {48}-3a_1\right) \frac 1 {n^4} +\left(-\frac {91}{2160}+a_1(6+4b_1)\right) \frac 1{n^5}+O(\frac 1{n^6}),
\end{align}
By Lemma 1, we know that the fastest possible sequence $(u_1(n))_{n\geq1}$ is obtained as the first item on the right of (2.8)
 vanishes. So taking $a_1=\frac 1 {144}, b_1=\frac 1{60}$, we can get the rate of convergence of the sequence $(u_1(n))_{n\geq1}$ is at least $n^{-5}$.

\noindent{\bf (Step 3) The second-correction.} We define the sequence $(u_2(n))_{n\geq1}$ by the relation
\begin{align}
W(n)=\sqrt \frac e \pi \left(1-\frac 1{2\left(n+\frac 13\right)}\right)^{n+\frac 13}\frac 1{\sqrt {n}}\exp \left(\frac 1{n^2}\frac {a_1}{n+b_1+\frac {a_2}{n+b_2}}\right)\exp u_2(n).
\end{align}

Using the same method as above, we obtain that the sequence $(u_2(n))_{n\geq1}$ converges fasted only if $a_2=\frac {781}{3600}, b_2=-\frac{4309}{109340}$.

\noindent{\bf (Step 4) The third-correction.} Similarly, define the sequence $(u_3(n))_{n\geq1}$ by the relation
\begin{align}
W(n)=\sqrt \frac e \pi \left(1-\frac 1{2\left(n+\frac 13\right)}\right)^{n+\frac 13}\frac 1{\sqrt {n}}\exp \left(\frac 1{n^2}\frac {a_1}{n+b_1+\frac {a_2}{n+b_2+\frac {a_3}{n+b_3}}}\right)\exp u_3(n)
\end{align}
Using the same method as above, we obtain that the sequence $(u_3(n))_{n\geq1}$ converges fasted only if $a_3=\frac {51396085}{89664267}, b3=\frac {25682346121}{449571834712}$.

The new asymptotic (1.5) is obtained.

\section{Proof of Theorem 2}
The double-side inequality (1.6) may be written as follows:
$$
f(n)=\ln W(n)-\frac 12+\frac 12\ln \pi-\left(n+\frac 13\right)\ln\left(1-\frac 1{2\left(n+\frac 13\right)}\right)-\frac 12\ln\frac 1n -\frac 1{n^2}\frac {\frac 1 {144}}{n+\frac 1{60}+\frac {\frac {781}{3600}}{n-\frac{4309}{109340}}}
$$
and
\begin{align*}
g(n)=&\ln W(n)-\frac 12+\frac 12\ln \pi-\left(n+\frac 13\right)\ln\left(1-\frac 1{2\left(n+\frac 13\right)}\right)\\
&-\frac 12\ln\frac 1n -\frac 1{n^2}\frac {\frac 1 {144}}{n+\frac 1{60}+\frac {\frac {781}{3600}}{n-\frac{4309}{109340}+\frac {\frac {51396085}{89664267}}{n+\frac {25682346121}{449571834712}}}}.
\end{align*}

Suppose $F(n)=f(n+1)-f(n)$ and $G(n)=g(n+1)-g(n)$. For every $x>1$, we can get
\begin{align}
F''(x)=\frac {A(x-1)}{160 x^4 (1 + x)^4 (1 + 2 x)^2 (1 + 3 x) (4 + 3 x) (6 x-1)^2 (5 +
   6 x)^2 \Psi_1^3(x;2) \Psi_2^3(x;2)}<0
\end{align}
and
\begin{align}
G''(x)=\frac {B(x)}{32 x^4 (1 + x)^4 (1 + 2 x)^2 (1 + 3 x) (4 + 3 x) (-1 + 6 x)^2 (5 +
   6 x)^2 \Psi_1^3(x;3) \Psi_2^3(x;3)}>0,
\end{align}
where
\begin{align*}
\Psi_1(x;2)=&10642 - 1119 x + 49203x^2,\\
\Psi_2(x;2)=&58726 + 97287 x + 49203 x^2, \\
\Psi_1(x;3)=&113505180 + 4083418255 x + 178132605 x^2 +   5180725368 x^3\\
\Psi_2(x;3)=&9555781408 + 19981859569 x + 15720308709 x^2 +   5180725368 x^3,
\end{align*}
$A(x)=-461003,...,512 n^{18}-...$ is a polynomial of 18th degree with all negative coefficients and $B(x)=212876,...,696x^{22}+...$ is a polynomial of 22th degree with all positive coefficients.

It show that $F(x)$ is strictly concave and $G(x)$ is strictly convex on $(0,\infty)$. According
to Theorem 1, when $n\rightarrow \infty$, it holds that $\lim_{n\rightarrow \infty}f (n) = \lim_{n\rightarrow \infty}g(n) = 0$; thus $\lim_{n\rightarrow \infty}F(n) = \lim_{n\rightarrow \infty}G(n) = 0$. As a result, we can make sure that $F(x) < 0$ and $G(x) > 0$ on $(0,\infty)$. Consequently, the sequence $f (n)$ is strictly increasing and
$g(n)$ is strictly decreasing while they both converge to $0$. As a result, we conclude that
$f (n) > 0$, and $g(n) < 0$ for every integer $n > 1$.

The proof of Theorem 2 is complete.

\section{Numerical computations}
In this section, we give Table 1 to demonstrate the superiority of our new series respectively. From what has been discussed above, we found out the new asymptotic function as follows:
\begin{align}
W(n)\sim \sqrt \frac e \pi \left(1-\frac 1{2\left(n+\frac 13\right)}\right)^{n+\frac 13}\frac 1{\sqrt {n}}\exp \left(\frac 1{n^2}\frac {\frac 1 {144}}{n+\frac 1{60}+\frac {\frac {781}{3600}}{n-\frac{4309}{109340}+\frac {\frac {51396085}{89664267}}{n+\frac {25682346121}{449571834712}}}}\right)=\sigma(n).
\end{align}

Chen and Qi \cite {chen2} gave:
\begin{align}
W(n)\sim \frac 1{\sqrt {\pi\left(n+\frac 14\right)}}=\alpha(n).
\end{align}

Qi and Mortici \cite {qi3} gave the improved formula:
\begin{align}
W(n)\sim \sqrt \frac e\pi \left(1-\frac 1{2n}\right)^n \sqrt \frac 1n\exp\left(\frac 1{24n^2}+\frac 1{48n^3}+\frac 1{160n^4}+\frac 1{960n^5}\right)=\beta(n)
\end{align}
and
\begin{align}
W(n)\sim \sqrt \frac e\pi \left(1-\frac 1{2\left(n+\frac 13\right)}\right)^{n+\frac 13} \sqrt \frac 1n=\gamma(n).
\end{align}

We can easily observe that the new formula converges fastest of the other three formulas.

\begin{table}
\centering
\caption{Simulations for $\alpha(n), \beta(n), \gamma(n)$ and $\sigma(n)$}
\begin{tabular}
{l c c c c}
\hline
n & $\frac{W(n)-\alpha(n)}{W(n)}$ & $\frac{W(n)-\beta(n)}{W(n)}$ & $\frac{W(n)-\gamma(n)}{W(n)}$ & $\frac{W(n)-\sigma(n)}{W(n)}$ \\
\hline
50   & $-6.1876\times 10^{-6}$ & $7.3576\times10^{-14}$  & $5.5532\times10^{-8}$ & $-3.8082\times10^{-19}$  \\
500  & $-6.2438\times 10^{-8}$ & $7.1643\times10^{-20}$  & $5.5554\times10^{-11}$ & $-3.8138\times10^{-28}$ \\
1000 & $-1.5617\times 10^{-8}$ & $1.1177\times10^{-21}$  & $6.9443\times10^{-12}$ & $-7.4489\times10^{-31}$ \\
2000 & $-3.9053\times 10^{-9}$ & $1.7452\times10^{-23}$  & $8.6805\times10^{-13}$ & $-1.4549\times10^{-33}$  \\
\hline
\end{tabular}
\end{table}

\section{Acknowledgements}
This work was supported by the National Natural Science Foundation of China (Grant No. 61403034), Beijing Municipal Commission of Education Science and Technology Program (KM201510017002).

\begin{flushleft}

Xu You \\
1. Department of Mathematics and Physics, Beijing Institute of Petrochemical Technology,\\
Beijing 102617, P. R. China \\
e-mail: youxu@bipt.edu.cn

\end{flushleft}
\end{document}